\documentclass{amsart}
\usepackage{amsmath}

\usepackage{mathrsfs}
\usepackage{amsfonts}
\usepackage{mathrsfs}
\usepackage{amssymb}
\usepackage{amsmath,amssymb}

\newtheorem{theorem}{Theorem}[section]
\newtheorem{lemma}[theorem]{Lemma}

\theoremstyle{definition}
\newtheorem{definition}[theorem]{Definition}

\theoremstyle{remark}
\newtheorem{remark}[theorem]{Remark}

\numberwithin{equation}{section}

\begin{document}

\title[ Riesz transform on Q-type space]
{Riesz transforms on Q-type spaces with application to
quasi-geostrophic equation }

\author{Pengtao\ Li}
\address{School of Mathematics, Peking University, Beijing, 100871, China}

\email{li\_ptao@163.com}

\author{Zhichun\ Zhai}
\address{Department of Mathematics and Statistics, Memorial University of Newfoundland, St. John's, NL A1C 5S7, Canada}
\curraddr{}
 \email{a64zz@mun.ca}
\thanks{Project supported in part  by Natural Science and
Engineering Research Council of Canada.}


\date{}

\dedicatory{}


\begin{abstract}
In this paper, we prove the boundedness of  Riesz transforms
$\partial_{j}(-\Delta)^{-1/2}$ ($j=1,2,\cdots,n$)  on the Q-type
spaces $Q_{\alpha}^{\beta}(\mathbb{R}^{n})$. As an application, we
get the well-posedness and regularity of the quasi-geostrophic
equation with initial data in $Q_{\alpha}^{\beta,
-1}(\mathbb{R}^{2})$.
\end{abstract}

\maketitle

\section{Introduction}
In this paper, we consider the boundedness of Riesz transforms on
the space $Q^{\beta}_{\alpha}(\mathbb{R}^{n})$, which was introduced
in \cite{Li Zhai} and defined as the set of all measurable functions
with
$$\sup_{I} (l(I))^{2\alpha-n+2\beta-2}\int_{I}\int_{I}\frac{|f(x)-f(y)|^{2}}{|x-y|^{n+2\alpha-2\beta+2}}dxdy<\infty$$
with the supremum being taken over all cubes $I$ with the edge
length $l(I)$ and the edges parallel to the coordinate axes in
$\mathbb{R}^{n}$. For $\beta=1$, the space
$Q^{\beta}_{\alpha}(\mathbb{R}^{n})$ becomes the classical space
$Q_{\alpha}(\mathbb{R}^{n})$  defined by the following norm:
\begin{equation}\label{eq2}
\|f\|_{Q_{\alpha}}=\sup_{I}\left((l(I))^{2\alpha-n}\int_{I}\int_{I}\frac{|f(x)-f(y)|^{2}}{|x-y|^{n+2\alpha}}dxdy\right)^{1/2}<\infty.
\end{equation}
This space was first introduced by M. Ess$\acute{e}$n, S. Janson, L.
Peng and J. Xiao in \cite{M. Essen S. Janson L. Peng J. Xiao}. As a
new space between Sobolev spaces $W^{1,n}(\mathbb{R}^{n})$ and
$BMO(\mathbb{R}^{n})$, it has been studied extensively by many
authors. We refer the readers to \cite{M. Essen S. Janson L. Peng J.
Xiao}, \cite{J. Xiao 1} and \cite{G. Dafni J. Xiao} for further
information and details.

Since the space $Q_{\alpha}(\mathbb{R}^{n})$ own a structure similar
to $BMO(\mathbb{R}^{n})$, it can be regarded as an analogy of
$BMO(\mathbb{R}^{n})$ in many cases. It is well-known that by the
equivalent characterization of Hardy space $H^{1}(\mathbb{R}^{n})$,
Riesz transforms $R_{j}=\partial_{j}(-\Delta)^{-1/2}$,
$j=1,2,\dots,n$ are bounded on $H^{1}(\mathbb{R}^{n})$. Then  the
duality between $H^{1}(\mathbb{R}^{n})$ and $BMO(\mathbb{R}^{n})$
obviously implies  the boundedness of
$R_{j}=\partial_{j}(-\Delta)^{-1/2}$ on $BMO(\mathbb{R}^{n})$. So it
is natural to ask if $R_{j},j=1,2,\dots,n$ are bounded on
$Q^{\beta}_{\alpha}(\mathbb{R}^{n})$. In Section 2, by an equivalent
characterization of $Q_{\alpha}^{\beta}(\mathbb{R}^{n})$ associated
with fractional heat semigroup $e^{-t(-\Delta)^{\beta}}$ obtained in
\cite{Li Zhai}, we prove that Riesz transforms $R_{j}$ are bounded
on the space $Q^{\beta}_{\alpha}(\mathbb{R}^{n})$. As far as we
know, our result is new even in the case
$Q_{\alpha}(\mathbb{R}^{n}),\alpha\in(0,1)$.

As an application, we consider the well-posedness  and regularity of
the quasi-geostrophic equation with initial data in
$Q^{\beta,-1}_{\alpha}(\mathbb{R}^{n})$. In recent years, Q-type
spaces have been applied to the study of PDE and Harmonic analysis
by several authors. For example, in \cite{J. Xiao 1}, J. Xiao
replaced $BMO^{-1}(\mathbb{R}^{n})$ in \cite{H. Koch D. Tataru} by a
new critical space $Q^{-1}_{\alpha}(\mathbb{R}^{n})$ which is
derivatives of $Q_{\alpha},\alpha\in(0,1)$ and got the
well-posedness of Naiver-Stokes equations with initial data in
$Q^{-1}_{\alpha}(\mathbb{R}^{n})$. When $\alpha=0$,
$Q^{-1}_{\alpha}(\mathbb{R}^{n})=BMO^{-1}(\mathbb{R}^{n})$, his
result generalized the well-posedness obtained by Koch and Tataru in
\cite{H. Koch D. Tataru}.

In \cite{Li Zhai}, inspiring by \cite{J. Xiao 1} and the scaling
invariant, we introduced a new Q-type space
$Q^{\beta}_{\alpha}(\mathbb{R}^{n})$ with $\alpha>0$,
$\max\{\frac{1}{2},\alpha\}<\beta<1$ such that $\alpha+\beta-1\geq0$
and considered the generalized Naiver-Stokes equations as follows.
\begin{equation}\label{eq3}
\left\{\begin{array} {l@{\quad \quad}l}
\partial_{t}u+(-\triangle)^{\beta}u+(u\cdot \nabla)u-\nabla p=0,
& \hbox{in}\ \mathbb{R}^{1+n}_{+};\\
\nabla \cdot u=0, & \hbox{in}\ \mathbb{R}^{1+n}_{+};\\
u|_{t=0}=u_{0}, &\hbox{in}\ \mathbb{R}^{n}. \end{array} \right.
\end{equation}
We proved  the well-posedness and regularity of the generalized
Naiver-Stokes equations with initial data in the space $Q^{\beta,\
-1}_{\alpha}$. For $\beta=1$, our spaces $Q^{\beta,\ -1}_{\alpha}$
retreat to $Q^{-1}_{\alpha}$ in \cite{J. Xiao 1}. So our result can
be regarded as a generalization of those of \cite{H. Koch D. Tataru}
and \cite{J. Xiao 1}.

In Section 3, We consider the two-dimensional subcritical
quasi-geostrophic
 dissipative equation $(DQG)_{\beta}$.
 \begin{equation}\label{eq1}
\left\{\begin{array} {l@{\quad \quad}l}
\partial_{t}\theta+(-\triangle)^{\beta}u+(u\cdot \nabla)\theta=0,
& \hbox{in}\ \mathbb{R}^{2}\times\mathbb{R}_{+}, \alpha>0;\\
u=\nabla^{\bot}(-\Delta)^{-1/2}\theta;\\
\theta(0,x)=\theta_{0}, &\hbox{in}\ \mathbb{R}^{2}. \end{array}
\right.
\end{equation}
where $\beta\in(\frac{1}{2},1)$, the scalar $\theta$ represent the
potential temperature, and $u$ is the fluid velocity.

The equations $(DQG)_{\beta}$ in either inviscid or dissipative
form, are special cases of the general quasi-geostrophic
approximations for the atmosphere and ocean flow with small Rossby
and Ekman numbers. Therefore, they are important models in
geophysical fluid dynamics. It was proposed by P. Constantin and A.
Majda, etc that the equations $(DQG)_{\beta}$ can be regarded as low
dimensional model equations for mathematical study of possible
development of singularity in smooth solutions of unforced
incompressible three dimensional fluid equations. See e.g.
\cite{CMT}, \cite{HPGS}, \cite{J1}, \cite{OY}, \cite{P} and the
references therein.

Recently, the equations $(DQG)_{\beta}$ have been intensively
studied because of their importance in mathematical and geophysical
fluid dynamics as mentioned above. Some important progress has been
made. We refer the readers to \cite{CH1}, \cite{CHL}, \cite{CMZ},
\cite{C}, \cite{CCW},  \cite{CW}, \cite{J2}, \cite{J Wu 1}, \cite{J
Wu 2} etc. for details.

In \cite{Marchand Rieusset}, F. Marchand and P. G.
Lemari$\acute{e}$-Rieusset studied the equations $(DQG)_{\beta}$ and
get the well-posedness of the solutions to the equation $(DQG)_{1}$
with the initial data in $BMO^{-1}(\mathbb{R}^{2})$. However,
because the space $BMO^{-1}(\mathbb{R}^{2})$ is invariant under the
scaling: $u_{0,\lambda}(x)=\lambda u_{0}(\lambda x)$, we see that
for the scaling corresponding to general $\beta<1$,
\begin{equation}\label{eq8}
\theta_{\lambda}(t,x)=\lambda^{2\beta-1}\theta(\lambda^{2\beta}t,
\lambda x), \qquad
\theta_{0,\lambda}(x)=\lambda^{2\beta-1}\theta_{0}(\lambda
x),\end{equation}
 the space $BMO^{-1}(\mathbb{R}^{2})$ is not invariant under this
scaling.

The above observation implies that if we want to generalize the
result in \cite{Marchand Rieusset} to general $\beta<1$, we should
choose a new space $X^{\beta}$ which satisfies the following two
properties. At first, the space $X^{\beta}$ should be invariant
under the scaling (\ref{eq8}). Secondly, $BMO^{-1}(\mathbb{R}^{2})$
is a ``special" case of $X^{\beta}$, that is, when $\beta=1$,
$X^{\beta}=BMO^{-1}(\mathbb{R}^{2})$.

In \cite{Li Zhai}, we have proved the space $Q^{\beta,\
-1}_{\alpha}(\mathbb{R}^{2})$ is exactly such a space. Therefore we
could apply our approach in \cite{Li Zhai} to the $(DQG)_{\beta}$
equation and get the well-posedness and regularity of the solution
to the $(DQG)_{\beta}$ equation.

It should be pointed out that the scope of $\beta$ is refined by the
choice of the space $Q^{\beta}_{\alpha}(\mathbb{R}^{n})$ In the
definition of $Q^{\beta}_{\alpha}(\mathbb{R}^{n}),$ the parameters
$\{\alpha,\ \beta\}$ should satisfy the condition: $\max\{\alpha,
\frac{1}{2}\}<\beta<1$ and $\alpha<\beta$ with $\alpha+\beta-1\geq0$
(see \cite{Li Zhai}). It is easy to see that $\beta>\frac{1}{2}$.

In \cite{RE}, the authors proved the global existence of the
solutions of the subcritical quasi-geostrophic equations with small
size initial data in the Besov norms paces
$\dot{B}^{1-2\beta,\infty}_{\infty}(\mathbb{R}^{2}).$ However our
well-posedness can't be deduced by the existence result in
\cite{RE}. In addition, owing to the structure of the
$Q^{\beta}_{\alpha}$, we can apply the method in \cite{Li Zhai} to
get the regularity of the solutions to the equation $(DQG)_{\beta}$.

 {\bf{Acknowledgements.}} We would like to thank our supervisor Professor Jie Xiao
  for discussion on this topic and kind  encouragement.

\section{Riesz transform on Q-type spaces $Q^{\beta}_{\alpha}(\mathbb{R}^{n})$}
In this section, we will prove that Riesz transforms are bounded on
Q-type spaces $Q^{\beta}_{\alpha}(\mathbb{R}^{n})$. At first we
recall the definition of $Q^{\beta}_{\alpha}(\mathbb{R}^{n})$.
\begin{definition}\label{def1}
Let $-\infty<\alpha$ and $\max\{\alpha,1/2\}<\beta<1$. Then $f\in
Q^{\beta}_{\alpha}(\mathbb{R}^{n})$ if and only if
$$\sup_{I}(l(I))^{2\alpha-n+2\beta-2}\int_{I}\int_{I}\frac{|f(x)-f(y)|^{2}}{|x-y|^{n+2\alpha-2\beta+2}}dxdy<\infty$$
where the supremum is taken over all cubes $I$ with the edge length
$l(I)$ and the edges parallel to the coordinate axes in
$\mathbb{R}^{n}$.
\end{definition}
For $\beta=1$ and $\alpha>-\infty$, the above spaces become the
$Q_{\alpha}(\mathbb{R}^{n})$ which were introduced in \cite{M. Essen
S. Janson L. Peng J. Xiao} by M. Essen, S. Janson, L. Peng and J.
Xiao. In \cite{G. Dafni J. Xiao}, G. Dafni and J. Xiao further
studied the structure of this space and get an equivalent
characterization via the heat semigroup associated with $\Delta$. It
has been proved that when $\alpha\in(0,1)$,
$Q_{\alpha}(\mathbb{R}^{n})\hookrightarrow BMO(\mathbb{R}^{n})$ and
when $\infty<\alpha<0$,
$Q_{\alpha}(\mathbb{R}^{n})=BMO(\mathbb{R}^{n})$(See \cite{M. Essen
S. Janson L. Peng J. Xiao}). Recall the definition of Morrey space
$\mathcal{L}_{p,\lambda}(\mathbb{R}^{n})$:
\begin{equation}
\|f\|_{\mathcal{L}_{p,\lambda}}=\sup_{I}\left((l(I))^{-\lambda}\int_{I}|f(x)-f_{I}|^{p}dx\right)^{1/p}<\infty.
\end{equation}
We see that when $\lambda=n$,
$\mathcal{L}_{p,\lambda}(\mathbb{R}^{n})=BMO(\mathbb{R}^{n})$ by
John-Nirenberg inequality. It is natural to ask if there exists some
relation between $\mathcal{L}_{p,\lambda}(\mathbb{R}^{n})$ and
$Q_{\alpha}(\mathbb{R}^{n})$. In fact, in \cite{J. Xiao 1}, J. Xiao
proved that for $\alpha\in(0,1)$,
$Q_{\alpha}(\mathbb{R}^{n})=(-\Delta)^{-\frac{\alpha}{2}}\mathcal{L}_{2,n-2\alpha}(\mathbb{R}^{n})$.

Following Xiao's idea in \cite{J. Xiao 1}, we will prove for our
space $Q^{\beta}_{\alpha}(\mathbb{R}^{n})$, a similar result holds.
At first we proved an equivalent characterization of
$\mathcal{L}_{2,n-2\gamma}(\mathbb{R}^{n})$ via the  semigroup
$e^{-t(-\Delta)^{\beta}}$. Here $e^{-t(-\Delta)^{\beta}}$ denotes
the convolution operator defined by Fourier transform:
$$\widehat{e^{-t(-\Delta)^{\beta}}f}(\xi)=e^{-t|\xi|^{2\beta}}\widehat{f}(\xi).$$
\begin{lemma}\label{le1} For $\gamma\in (0,1)$.
Let f be a
measurable complex-valued function on $\mathbb{R}^{n}$. Then
$$f\in \mathcal{L}_{2, n-\gamma}(\mathbb{R}^{n})\Longleftrightarrow \sup_{x\in\mathbb{R}^{n}, r\in(0,\infty)}r^{2\gamma-n}\int_{0}^{r}\int_{|y-x|<r}
\left|\nabla
e^{-t^{2\beta}(-\beta)^{\beta}}f(y)\right|^{2}tdydt<\infty.$$
\end{lemma}
\begin{proof}
Taking $(\psi_{0})_{t}(x)=t\nabla
e^{-t^{2\beta}(-\Delta)^{\beta}}(x,0)$ with the Fourier symbol
$\widehat{(\psi_{0})_{t}(x)}(\xi)=t|\xi|e^{-t^{2\beta}|\xi|^{2\beta}}$.
Define  a ball $B=\{y\in\mathbb{R}^{n}:|y-x|<r\}$  and
$f_{2B}=\frac{1}{|2B|}\int_{2B}f(x)dx$ is the mean of $f$ on $2B$.
We split $f$ into $f=f_{1}+f_{2}+f_{3}$ where
$f_{1}=(f-f_{2B})\chi_{2B}$, $f_{2}=(f-f_{2B})\chi_{(2B)^{c}}$ and
$f_{3}=f_{2B}$. Because
$$\int_{}^{}(\psi_{0})_{t}(x)dx=\int t\nabla e^{-t^{2\beta}(-\Delta)^{\beta}}(x,0)dx=0,$$
we have
$$t\nabla e^{-t^{2\beta}(-\Delta)^{\beta}}f(y)=(\psi_{0})_{t}\ast f(y)=(\psi_{0})_{t}\ast f_{1}(y)+(\psi_{0})_{t}\ast f_{2}(y).$$
It is easy to see that
\begin{eqnarray*}
\int_{0}^{r}\int_{B}\left|(\psi_{0})_{t}\ast
f_{1}(y)\right|^{2}\frac{dydt}{t}&\lesssim&\int_{0}^{r}\int_{\mathbb{R}^{n}}\left|(\psi_{0})_{t}\ast
f_{1}(y)\right|^{2}\frac{dydt}{t}\\
&=&\left\|\left(\int_{0}^{\infty}\left|(\psi_{0})_{t}\ast
f_{1}(\cdot)\right|^{2}\frac{dt}{t}\right)^{1/2}\right\|_{L^{2}(dy)}.
\end{eqnarray*}
Because $(\psi_{0})_{1}=\nabla e^{-(-\Delta)^{\beta}}$, obviously we
have $\int (\psi_{0})_{1}(x)dx=1$ and $(\psi_{0})_{1}$ belongs to
the Schwartz class $\mathcal{S}$, the function
$$G(f)=\left(\int_{0}^{\infty}\left|(\psi_{0})_{t}\ast f_{1}(y)\right|^{2}\frac{dt}{t}\right)^{1/2}$$
is a Littlewood-Paley G-function. So we can get
\begin{eqnarray*}
\int_{0}^{r}\int_{B}\left|(\psi_{0})_{t}\ast
f_{1}(y)\right|^{2}\frac{dydt}{t}&\lesssim&\int_{2B}\left|f(y)-f_{2B}\right|^{2}dy\\
&\lesssim&r^{n-2\gamma}\|f\|^{2}_{\mathcal{L}_{2,n-2\gamma}(\mathbb{R}^{n})}.
\end{eqnarray*}
Now we estimate the term associated with $f_{2}(y)$. Because
\begin{eqnarray*}
|(\psi_{0})_{t}\ast f_{2}(y)|&=&\left|\int_{\mathbb{R}^{n}}t\nabla
e^{-t^{2\beta}(-\Delta)^{\beta}}(y-z)f_{2}(z)dz\right|\\
&\lesssim&\int_{\mathbb{R}^{n}\backslash 2B}|t\nabla
e^{-t^{2\beta}(-\Delta)^{\beta}}(y-z)||f(z)-f_{2B}|dz\\
&\lesssim&\int_{\mathbb{R}^{n}\backslash
2B}\frac{t|f(z)-f_{2B}|}{t^{n+1}(1+t^{-1}|z-y|)^{n+1}}dz
\end{eqnarray*}
where in the last inequality, we have used the following estimate:
$$\left|\nabla e^{-t(-\Delta)^{\beta}}(x, y)\right|\lesssim\frac{1}{t^{\frac{n+1}{2\beta}}}\frac{1}{(1+t^{-\frac{1}{2\beta}}|x-y|)^{n+1}}.$$
Set $B_{k}=B(x,2^{k})$. For every $(t,y)\in(0,r)\times B(x,r)$, we
have $0<t<r$ and $|x-y|<r$. If $z\in B_{k+1}\backslash B_{k}$, that
is, $|z-x|>2^{k}r$, we have $|x-y|<|x-z|/2$ and
\begin{eqnarray*}
|(\psi_{0})_{t}\ast
f_{2}(y)|&\lesssim&\int_{\mathbb{R}^{n}\backslash
2B}\frac{t|f(z)-f_{2B}|}{(t+|z-x|)^{n+1}}dz\\
&\lesssim&\sum_{k=1}^{\infty}\frac{t}{(2^{k}r)^{n+1}}\int_{2^{k+1}B}|f(z)-f_{2B}|dz\\
&\lesssim&
t\sum_{k=1}^{\infty}\frac{(2^{k+1}r)^{n}}{(2^{k}r)^{n+1}}\left(\frac{1}{(2^{k+1}r)^{n}}\int_{2^{k+1}B}|f(z)-f_{2B}|^{2}dz\right)^{1/2}\\
&\lesssim&t\sum_{k=1}^{\infty}\frac{1}{2^{k}r}\left[\left(\frac{1}{(2^{k+1}r)^{n}}\int_{2^{k+1}B}|f(z)-f_{2^{k+1}B}|^{2}dz\right)^{1/2}+\left|f_{2^{k+1}B}-f_{2B}\right|\right]\\
&\lesssim&t\left[\sum_{k=1}^{\infty}\frac{1}{2^{k}r}\left(\frac{1}{(2^{k+1}r)^{n}}\int_{2^{k+1}B}|f(z)-f_{2^{k+1}B}|^{2}dz\right)^{1/2}
+\sum_{k=1}^{\infty}\frac{1}{2^{k}r}\left|f_{2^{k+1}B}-f_{2B}\right|\right]\\
&=:&t(S_{1}+S_{2}).
\end{eqnarray*}
For $S_{1}$, we have
\begin{eqnarray*}
S_{1}&=&t\sum_{k=1}^{\infty}\frac{1}{2^{k}r}
\left(\frac{(2^{k+1}r)^{n-2\gamma}}{(2^{k+1}r)^{n}}\frac{1}{(2^{k+1}r)^{n-2\gamma}}\int_{2^{k+1}B}|f(z)-f_{2^{k+1}B}|^{2}dz\right)^{1/2}\\
&\lesssim&t\sum_{k=1}^{\infty}\frac{1}{2^{k}r}r^{-\gamma}\|f\|_{\mathcal{L}_{2,n-2\gamma}(\mathbb{R}^{n})}\\
&\lesssim&tr^{-1-\gamma}\|f\|_{\mathcal{L}_{2,n-2\gamma}(\mathbb{R}^{n})}.
\end{eqnarray*}
For $S_{2}$, we have
$$S_{2}\lesssim t\sum_{k=1}^{\infty}\frac{1}{2^{k}r}[|f_{2B}-f_{4B}|+\cdots+|f_{2^{k}B}-f_{2^{k+1}B}|].$$
For $\forall 2\leq j\leq k$, it is easy to see that
\begin{eqnarray*}
|f_{2^{j}B}-f_{2^{j+1}B}|&\lesssim&\frac{1}{|2^{j}B|}\int_{2^{j}B}|f(z)-f_{2^{j+1}B}|dz\\
&\lesssim&\left(\frac{1}{|2^{j}B|}\int_{2^{j}B}|f(z)-f_{2^{j+1}B}|^{2}dz\right)^{1/2}\\
&\lesssim&r^{-\gamma}\|f\|_{\mathcal{L}_{2,n-2\gamma}(\mathbb{R}^{n})}.
\end{eqnarray*}
Then we can get
\begin{eqnarray*}
S_{2}\lesssim t\sum_{k=1}^{\infty}\frac{1}{2^{k}r}k\cdot
r^{-\gamma}\|f\|_{\mathcal{L}_{2,n-2\gamma}(\mathbb{R}^{n})}\lesssim
tr^{-1-\gamma}\|f\|_{\mathcal{L}_{2,n-2\gamma}(\mathbb{R}^{n})}.
\end{eqnarray*}
Therefore we can get
\begin{eqnarray*}
\int_{0}^{r}\int_{B}\left|(\psi_{0})_{t}\ast
f_{2}(y)\right|^{2}t^{-1}dydt &\lesssim&
\int_{0}^{r}\int_{B}t^{2}r^{-2\gamma-2}\|f\|^{2}_{\mathcal{L}_{2,n-2\gamma}(\mathbb{R}^{n})}dydt\\
&\lesssim&\|f\|^{2}_{\mathcal{L}_{2,n-2\gamma}(\mathbb{R}^{n})}r^{-2\gamma-2}|B|\int_{0}^{r}tdt\\
&\lesssim
&r^{n-2\gamma}\|f\|^{2}_{\mathcal{L}_{2,n-2\gamma}(\mathbb{R}^{n})}.
\end{eqnarray*}
For the converse, let $S(I)=\{(t,x)\in\mathbb{R}^{n+1}_{+},
0<t<l(I), x\in I\}$ if $f$ such that
\begin{eqnarray*}
&&\sup_{I}[l(I)]^{2\gamma-n}\int_{S(I)}\left|t\nabla
e^{-t^{2\beta}(-\Delta)^{\beta}}f(y)\right|^{2}\frac{dydt}{t}\\
&&=\sup_{I}[l(I)]^{2\gamma-n}\int_{S(I)}\left|\nabla
e^{-t^{2\beta}(-\Delta)^{\beta}}f(y)\right|^{2}tdydt<\infty.
\end{eqnarray*}
Denote
$$\Pi_{\psi_{0}}F(x)=\int_{\mathbb{R}^{n+1}_{+}}F(t,y)(\psi_{0})_{t}(x-y)\frac{dydt}{t},$$
we will prove that if
$$\|F\|_{C_{\gamma}}=\sup_{I}\left([l(I)]^{2\gamma-n}\int_{S(I)}\left|F(t,y)\right|^{2}\frac{dydt}{t}\right)^{1/2}<\infty,$$
then for any cube $J\subset\mathbb{R}^{n}$
$$\int_{J}\left|\Pi_{\psi_{0}}F(x)-(\Pi_{\psi_{0}}F)_{J}\right|^{2}dx\lesssim [l(J)]^{n-2\gamma}\|F\|_{C_{\gamma}}^{2}.$$
We split $F$ into
$F=F_{1}+F_{2}=F\mid_{S(2J)}+F\mid_{\mathbb{R}^{n+1}\setminus
S(2J)}$ and get
\begin{eqnarray*}
\int_{J}\left|\Pi_{\psi_{0}}F_{1}(x)\right|^{2}dx&\leq&\int_{J}\left|\Pi_{\psi_{0}}F_{1}(x)\right|^{2}dx\\
&\leq&\int_{S(2J)}\left|F(t,y)\right|^{2}\frac{dydt}{t}\\
&\lesssim&[l(J)]^{n-2\gamma}\|F\|^{2}_{C_{\gamma}}.
\end{eqnarray*}
Now we estimate the term associated with $F_{2}$. We have
\begin{eqnarray*}
&&\int_{J}\left|\Pi_{\psi_{0}}F_{1}(x)\right|^{2}dx\\
&&=\int_{J}\left|\int_{\mathbb{R}^{n+1}_{+}}(\psi_{0})_{t}(x-y)F_{2}(t,y)t^{-1}dydt\right|^{2}dx\\
&&\lesssim\int_{J}\left(\int_{\mathbb{R}^{n+1}_{+}\setminus
S(2J)}|(\psi_{0})_{t}(x-y)||F_{2}(t,y)|\frac{dydt}{t}\right)^{2}dx\\
&=&\int_{J}\left(\sum_{k=1}^{\infty}\int_{S(2^{k+1}J)\setminus
S(2^{k}J)}|(\psi_{0})_{t}(x-y)||F_{2}(t,y)|\frac{dydt}{t}\right)^{2}dx.
\end{eqnarray*}
Because
$$|(\psi_{0})_{t}(x-y)|\lesssim\frac{t}{t^{n+1}(1+t^{-1}|x-y|)^{n+1}},$$
we have
\begin{eqnarray*}
&&\int_{J}\left|\Pi_{\psi_{0}}F_{1}(x)\right|^{2}dx\\
&&\lesssim\int_{J}\left(\sum_{k=1}^{\infty}\int_{S(2^{k+1}J)\setminus
S(2^{k}J)}\frac{t}{[t+2^{k}l(J)]^{n+1}}|F_{2}(t,y)|\frac{dydt}{t}\right)^{2}dx\\
&&\lesssim\int_{J}\left(\sum_{k=1}^{\infty}(2^{k}l(J))^{-(n+1)}\int_{S(2^{k+1}J)\setminus
S(2^{k}J)}|F_{2}(t,y)|dydt\right)^{2}dx\\
&&\lesssim\int_{J}\left[\sum_{k=1}^{\infty}[2^{k}l(J)]^{-n}[2^{k+1}l(J)]^{n/2}\left(\int_{S(2^{k+1}J)\setminus
S(2^{k}J)}|F_{2}(t,y)|^{2}\frac{dydt}{t}\right)^{1/2}\right]^{2}dx\\
&&\lesssim\|F\|^{2}_{C_{\gamma}}[l(J)]^{n-2\gamma}.
\end{eqnarray*}
Therefore we get
\begin{eqnarray*}
\int_{J}\left|\Pi_{\psi_{0}}F(x)-(\Pi_{\psi_{0}}F)_{J}\right|^{2}dx&\leq&\int_{J}\left|\Pi_{\psi_{0}}F(x)\right|^{2}dx\\
&\lesssim&\int_{J}\left|\Pi_{\psi_{0}}F_{1}(x)\right|^{2}dx+\int_{J}\left|\Pi_{\psi_{0}}F_{2}(x)\right|^{2}dx\\
&\lesssim&\|F\|^{2}_{C_{\gamma}}[l(J)]^{n-2\gamma}.
\end{eqnarray*}
Because
$$\Pi_{\psi_{0}}F(x)=\int(\psi_{0})_{t}\ast(\psi_{0})_{t}\ast f\frac{dt}{t},$$
by Calder$\acute{o}$n reproducing formula, we have
$\Pi_{\psi_{0}}F(x)=f(x)$, that is, $f(x)=\Pi_{\psi_{0}}F(x)\in
\mathcal{L}_{2,n-2\gamma}(\mathbb{R}^{n})$. This completes the proof
of Lemma \ref{le1}.
\end{proof}

\begin{theorem}\label{th1}
For $\alpha>0$, $\max\{\alpha,\frac{1}{2}\}<\beta<1$ with
$\alpha+\beta-1\geq0$, we have
$$Q^{\beta}_{\alpha}(\mathbb{R}^{n})=(-\Delta)^{-\frac{(\alpha-\beta+1)}{2}}\mathcal{L}_{2,\
n-2(\alpha+\beta-1)}(\mathbb{R}^{n}).$$
\end{theorem}
\begin{proof}
For $f\in \mathcal{L}_{2,\ n-2(\alpha+\beta-1)}(\mathbb{R}^{n}).$
Let $F(t,y)=t^{\alpha-\beta+1}t\nabla
e^{-t^{2\beta}(-\Delta)^{\beta}}f(y)$. By Lemma \ref{le1}, we have
\begin{eqnarray*}
&&r^{2(\alpha+\beta-1)-n}\int_{0}^{r}\int_{|y-x|<r}|F(t,y)|^{2}\frac{dydt}{t^{1+2(\alpha-\beta+1)}}\\
&&\lesssim
r^{2(\alpha+\beta-1)-n}\int_{0}^{r}\int_{|y-x|<r}|t^{\alpha-\beta+1}t\nabla
e^{-t^{2\beta}(-\Delta)^{\beta}}f(y)|^{2}\frac{dydt}{t^{1+2(\alpha-\beta+1)}}\\
&&\lesssim r^{2(\alpha+\beta-1)-n}\int_{0}^{r}\int_{|y-x|<r}|t\nabla
e^{-t^{2\beta}(-\Delta)^{\beta}}f(y)|^{2}\frac{dydt}{t}\\
&&\lesssim\|f\|_{\mathcal{L}_{2,\
n-2(\alpha+\beta-1)}(\mathbb{R}^{n})}.
\end{eqnarray*}
So $F\in T^{\infty}_{\alpha,\beta}$ where the space
$T^{\infty}_{\alpha,\beta}$ is a tent space defined in \cite{Li
Zhai}(See \cite{Li Zhai}, Definition 3.5 for details). By Theorem
3.21 in \cite{Li Zhai}, $\Pi_{\psi_{0}}$ is bounded from
$T^{\infty}_{\alpha,\beta}$ to $Q^{\beta}_{\alpha}$. Therefore we
have
$$\|f\|_{Q^{\beta}_{\alpha}(\mathbb{R}^{n})}=\|\Pi_{\psi_{0}}F\|_{Q^{\beta}_{\alpha}(\mathbb{R}^{n})}\lesssim \|F\|_{T^{\infty}_{\alpha,\beta}}.$$
We have
$$\widehat{F}(t,\ \xi)=t^{\alpha-\beta+2}|\xi|e^{-t^{2\beta}|\xi|^{2\beta}}\widehat{f}(\xi),$$
then
\begin{eqnarray*}
\widehat{\Pi_{\psi_{0}}F}(\xi)&=&\int_{0}^{\infty}\widehat{F}(t,\xi)\widehat{(\psi_{0})_{t}}(\xi)\frac{dt}{t}\\
&=&\int_{0}^{\infty}t^{\alpha-\beta+2}|\xi|e^{-t^{2\beta}|\xi|^{2\beta}}t|\xi|e^{-t^{2\beta}|\xi|^{2\beta}}\widehat{f}(\xi)\frac{dt}{t}\\
&=&|\xi|^{2}\widehat{f}(\xi)\int_{0}^{\infty}t^{\alpha-\beta+2}e^{-t^{2\beta}|\xi|^{2\beta}}dt.
\end{eqnarray*}
Set $t^{2\beta}=s$ and $|\xi|^{2\beta}s=u$. We can get
\begin{eqnarray*}
\widehat{\Pi_{\psi_{0}}F}(\xi)
&=&\int^{\infty}_{0}s^{\frac{\alpha-\beta+2}{2\beta}}e^{-2s|\xi|^{2\beta}}s^{\frac{1}{2\beta}-1}ds\widehat{f}(\xi)|\xi|^{2}\\
&=&\widehat{f}(\xi)|\xi|^{2}\int_{0}^{\infty}(u|\xi|^{-2\beta})^{\frac{\alpha-\beta+3}{2\beta}-1}e^{-u}|\xi|^{-2\beta}du\\
&=&\widehat{f}(\xi)|\xi|^{2}|\xi|^{-(\alpha-\beta+3)+2\beta-2\beta}\int_{0}^{\infty}u^{\frac{\alpha-\beta+3}{2\beta}-1}e^{-2u}du.
\end{eqnarray*}
Because $\frac{1}{2}<\beta<1$ and $0<\alpha<\beta$, the integral
$\int^{\infty}_{0}u^{\frac{\alpha-\beta+3}{2\beta}-1}e^{-2u}du$ is a
constant. We denote it by $C_{\alpha,\beta}$, so
$$\widehat{\Pi_{\psi_{0}}F}(\xi)=C_{\alpha,\beta}\widehat{f}(\xi)|\xi|^{-(\alpha-\beta+1)}.$$
By inverse Fourier transform, we have
$$\Pi_{\psi_{0}}F(x)=C_{\alpha,\beta}(-\Delta)^{-\frac{\alpha-\beta+1}{2}}f(x).$$
Conversely, suppose $g\in Q^{\beta}_{\alpha}(\mathbb{R}^{n})$. Set
$G(t,y)=t^{1-(\alpha-\beta+1)}\nabla
e^{-t^{2\beta}(-\Delta)^{\beta}}g(y)$. We have, by the equivalent
characterization of $Q^{\beta}_{\alpha}(\mathbb{R}^{n})$ (see
\cite{Li Zhai} for details).
\begin{eqnarray*}
&&\left([l(I)]^{2(\alpha+\beta-1)-n}\int_{S(I)}\left|t^{1-2(\alpha-\beta+1)}\nabla
e^{-t^{2\beta}(-\Delta)^{\beta}}g(y)\right|^{2}\frac{dydt}{t}\right)^{1/2}\\
&&=\left([l(I)]^{2(\alpha+\beta-1)-n}\int_{S(I)}\left|t\nabla
e^{-t^{2\beta}(-\Delta)^{\beta}}g(y)\right|^{2}\frac{dydt}{t^{1+2(\alpha-\beta+1)}}\right)^{1/2}\\
&&\lesssim \|g\|_{Q^{\beta}_{\alpha}(\mathbb{R}^{n})},
\end{eqnarray*}
that is, $G(t,y)\in C_{\alpha+\beta-1}$. By Lemma \ref{le1}, we have
$\Pi_{\psi_{0}}G(t,y)\in\mathcal{L}_{2,\ n-2(\alpha+\beta-1)}.$ We
have
\begin{eqnarray*}
\widehat{f}(\xi)&=&\widehat{\Pi_{\psi_{0}}G}(t,\xi)\\
&=&\int_{0}^{\infty}t|\xi|e^{-t^{2\beta}|\xi|^{2\beta}}t^{1-(\alpha-\beta+1)}|\xi|e^{-t^{2\beta}|\xi|^{2\beta}}\widehat{g}(\xi)\frac{dt}{t}\\
&=&C_{\alpha,\beta}|\xi|^{1+(\alpha-\beta)}\widehat{g}(\xi)\\
&=&C_{\alpha,\beta}\widehat{[(-\Delta)^{\frac{\alpha-\beta+1}{2}}g]}(\xi).
\end{eqnarray*}
Then $f(x)=C_{\alpha,\beta}(-\Delta)^{\frac{\alpha-\beta+1}{2}}g$.
Thus,  we get
$Q^{\beta}_{\alpha}(\mathbb{R}^{n})=(-\Delta)^{-\frac{(\alpha-\beta+1)}{2}}\mathcal{L}_{2,\
n-2(\alpha+\beta-1)}(\mathbb{R}^{n}).$
\end{proof}
\begin{theorem}\label{th2}
Suppose $\alpha>0$, $\max{\alpha,\frac{1}{2}}<\beta<1$ with
$\alpha+\beta-1\geq0$. For $j=1,2,\dots , n$,
  the Riesz transforms $R_{j}=\partial_{j}(-\Delta)^{-1/2}$   are bounded on the $Q-$type spaces
$Q^{\beta}_{\alpha}(\mathbb{R}^{n})$.
\end{theorem}
\begin{proof}
Notice the equivalent norm of $Q^{\beta}_{\alpha}(\mathbb{R}^{n})$.
 $f\in Q^{\beta}_{\alpha}(\mathbb{R}^{n})$ if and only if
 \begin{equation}\label{char 1}
 \sup_{x\in\mathbb{R}^{n},r\in(0,\infty)}r^{2\alpha-n+2\beta-2}
\int_{0}^{r^{2\beta}}\int_{|y-x|<r}|\nabla
e^{-t(-\triangle)^{\beta}}
f(y)|^{2}t^{-\frac{\alpha}{\beta}}dydt<\infty.
\end{equation}
As  a convolution operator, Riesz transform $R_{j}$ and $\nabla$ can
change the order of operation. So we only need to estimate the term
$$r^{2\alpha-n+2\beta-2}\int_{0}^{r^{2\beta}}\int_{|y-x_{0}|<r}\left|R_{j}f(t,y)\right|^{2}\frac{dydt}{t^{\alpha/\beta}}$$
where $0<\alpha<\beta$, $\alpha+\beta-1\geq0$ and $j=1,2,\dots,n$.
We split $f(t,y)$ into
$$f(t,y)=f_{0}(t,y)+\sum_{k=1}^{\infty}f_{k}(t,y),$$
where $f_{0}(t,y)=f(t,y)\chi_{B(x_{0},2r)}(y)$ and
$f_{k}(t,y)=f(t,y)\chi_{B(x_{0},2^{k+1}r)\setminus
B(x_{0},2^{k}r)}(y)$. We have
\begin{eqnarray*}
&&\left(r^{2\alpha-n+2\beta-2}\int_{0}^{r^{2\beta}}\int_{|y-x_{0}|<r}\left|R_{j}f(t,y)\right|^{2}\frac{dydt}{t^{\alpha/\beta}}\right)^{1/2}\\
&&\leq\left(r^{2\alpha-n+2\beta-2}\int_{0}^{r^{2\beta}}\int_{|y-x_{0}|<r}\left|R_{j}f_{0}(t,y)\right|^{2}\frac{dydt}{t^{\alpha/\beta}}\right)^{1/2}\\
&&+\sum_{k=1}^{\infty}\left(r^{2\alpha-n+2\beta-2}\int_{0}^{r^{2\beta}}
\int_{|y-x_{0}|<r}\left|R_{j}f(t,y)\right|^{2}\frac{dydt}{t^{\alpha/\beta}}\right)^{1/2}\\
&&:=M_{0}+\sum_{k=1}^{\infty}M_{k}.
\end{eqnarray*}
By the $L^{2}-$boundedness of Riesz transforms $R_{j},\
j=1,2,\dots,n$, we have
\begin{eqnarray*}
M_{0}&\lesssim&\left(r^{2\alpha-n+2\beta-2}\int_{0}^{r^{2\beta}}\int_{|y-x_{0}|<r}\left|f(t,y)\right|^{2}\frac{dydt}{t^{\alpha/\beta}}\right)^{1/2}\\
&\lesssim&C_{\alpha}\sup_{x_{0}\in\mathbb{R}^{n},\ r>0}
\left(r^{2\alpha-n+2\beta-2}\int_{0}^{r^{2\beta}}\int_{|y-x_{0}|<r}\left|f(t,y)\right|^{2}\frac{dydt}{t^{\alpha/\beta}}\right)^{1/2}\\
\end{eqnarray*}
Now we estimate $M_{k}$. We only need to estimate the integration as
follows.
$$I=\int_{|y-x_{0}|<r}\left|R_{j}f_{k}(t,y)\right|^{2}dy.$$
As a singular integral operator,
$$R_{j}g(x)=\int_{\mathbb{R}^{n}}\frac{x_{j}-y_{j}}{|x-y|^{n+1}}g(y)dy.$$
By H\"{o}lder's inequality, we can get
\begin{eqnarray*}
I&=&\int_{|y-x_{0}|<r}\left|\int_{2^{k}r\leq|z-x_{0}|<2^{k+1}r}\frac{y_{j}-z_{j}}{|y-z|^{n+1}}f(t,z)dz\right|^{2}dy\\
&\lesssim&\int_{|y-x_{0}|<r}\left(\frac{1}{(2^{k}r)^{n}}\int_{|z-x_{0}|<2^{k+1}r}|f(t,z)|dz\right)^{2}dy\\
&\lesssim&\int_{|y-x_{0}|<r}\frac{1}{(2^{k}r)^{n}}\int_{|z-x_{0}|<2^{k+1}r}|f(t,z)|^{2}dzdy\\
&\lesssim&\frac{1}{2^{kn}}\int_{|z-x_{0}|<2^{k+1}r}|f(t,z)|^{2}dz.
\end{eqnarray*}
So we have
\begin{eqnarray*}
M_{k}&=&\left(r^{2\alpha-n+2\beta-2}\int_{0}^{r^{2\beta}}
\int_{|y-x_{0}|<r}\left|R_{j}f(t,y)\right|^{2}\frac{dydt}{t^{\alpha/\beta}}\right)^{1/2}\\
&\lesssim&\left(r^{2\alpha-n+2\beta-2}\frac{1}{2^{kn}}\int_{0}^{r^{2\beta}}
\int_{|z-x_{0}|<2^{k+1}r}\left|f(t,z)\right|^{2}\frac{dzdt}{t^{\alpha/\beta}}\right)^{1/2}\\
&\lesssim&\left(2^{-k(2\alpha-n+2\beta-2)}\frac{1}{2^{kn}}(2^{k}r)^{2\alpha-n+2\beta-2}\int_{0}^{r^{2\beta}}
\int_{|z-x_{0}|<2^{k+1}r}\left|f(t,z)\right|^{2}\frac{dzdt}{t^{\alpha/\beta}}\right)^{1/2}\\
&\lesssim&(2^{-k(2\alpha-n+2\beta-2)}\frac{1}{2^{kn}})^{1/2}\sup_{x_{0}\in\mathbb{R}^{n},\
r>0}\left(r^{2\alpha-n+2\beta-2}\int_{0}^{r^{2\beta}}\int_{|z-x_{0}|<r}|f(t,z)|^{2}\frac{dzdt}{t^{\alpha/\beta}}\right)^{1/2}.
\end{eqnarray*}
Therefore we can get
\begin{eqnarray*}
&&M_{0}+\sum_{k=1}^{\infty}M_{k}\\
&&\lesssim[1+\sum_{k=1}^{\infty}2^{-k(\alpha+\beta-1)}]\sup_{x_{0}\in\mathbb{R}^{n},\
r>0}\left(r^{2\alpha-n+2\beta-2}\int_{0}^{r^{2\beta}}\int_{|z-x_{0}|<r}|f(t,z)|^{2}\frac{dzdt}{t^{\alpha/\beta}}\right)^{1/2}\\
&&\lesssim C\|f\|_{Q^{\beta}_{\alpha}(\mathbb{R}^{n})}.
\end{eqnarray*}
This completes the proof of Theorem \ref{th2}.

\end{proof}
Similar to the proof of Theorem \ref{th2}, we can prove the
following theorem.
\begin{theorem}\label{th3}
For a singular operator $T$ defined by
$$Tf(x)=\int_{\mathbb{R}^{n}}K(x-y)f(y)dy$$
where the kernel $K(x)$ satisfies:
$$|\partial^{\gamma}_{x}K(x)|\leq A_{\gamma}|x|^{-n-\gamma},\quad (\gamma>0)$$
or equivalently $\widehat{Tf}(\xi)=m(\xi)\widehat{f}(\xi)$ with the
symbol $m(\xi)$ satisfies:
$$|\partial^{\gamma}_{\xi}m(\xi)|\leq A_{\gamma'}|\xi|^{-\gamma}$$
holds for all $\gamma$. Suppose $\alpha>0$,
$\max\{\alpha,\frac{1}{2}\}<\beta<1$ with $\alpha+\beta-1\geq0$, we
have $T$ is bounded on the $Q-$type spaces
$Q^{\beta}_{\alpha}(\mathbb{R}^{n})$.
\end{theorem}
\section{Well-posedness and regularity of quasi-geostrophic equation}
In this section, we study the well-posedness and regularity of
quasi-geostrophic equation with initial data in the space
$Q^{\beta}_{\alpha}(\mathbb{R}^{2})$. We introduce the definition of
$X^{\beta}_{\alpha}(\mathbb{R}^{n})$.
\begin{definition}\label{def2}
The space $X^{\beta}_{\alpha}(\mathbb{R}^{2})$ consists of the
functions which are locally integrable on
$(0,\infty)\times\mathbb{R}^{2}$ such that
$$\sup_{t>0}t^{1-\frac{1}{2\beta}}\|f(t,\cdot)\|_{\dot{B}^{0,1}_{\infty}(\mathbb{R}^{2})}<\infty$$
and
$$\sup_{x\in\mathbb{R}^{2},\ r>0}r^{2\alpha+2\beta-4}\int_{0}^{r^{2\beta}}\int_{|y-x_{0}|<r}
|f(t,y)|^{2}+|R_{1}f(t,y)|^{2}+|R_{2}f(t,y)|^{2}\frac{dydt}{t^{\alpha/\beta}}<\infty,$$
where $R_{j},\ (j=1,2)$ denote the Riesz transforms in
$\mathbb{R}^{2}$.
\end{definition}
For the quasi-geostrophic dissipative equations

\begin{equation} \label{eq4}
\left\{ \begin{aligned}
         \partial_{t}\theta &= -(-\Delta)^{\beta}+\partial_{1}(\theta R_{2}\theta)-\partial_{2}(\theta R_{1}\theta); \\
                  \theta(0,x)&=\theta_{0}(x)
                          \end{aligned} \right.
                          \end{equation}
where $\beta\in(\frac{1}{2},1)$. The solution to the equation
(\ref{eq4}) can be represented as
$$u(t,x)=e^{-t(-\Delta)^{\beta}}u_{0}+B(u,u)$$
where the bilinear form $B(u,v)$ is defined by
$$B(u,v)=\int_{0}^{t}e^{-(t-s)(-\Delta)^{\beta}}(\partial_{1}(vR_{2}u)-\partial_{2}(vR_{1}u))ds.$$
In order to prove the well-posedness, we need the following
preliminary lemmas. For their proof, we refer the readers to
\cite{Li Zhai}, Lemma 4.8 and Lemma 4.9.
\begin{lemma} \label{le 2} Given $\alpha\in(0,1).$ For a fixed $T\in (0,\infty]$
and a function $f(\cdot, \cdot)$ on $\mathbb{R}^{1+n}_{+},$ let
$A(t)=\int_{0}^{t}e^{-(t-s)(-\triangle)^{\beta}}(-\triangle)^{\beta}f(s,x)ds.$
Then
\begin{equation}\label{eq1c}
\int_{0}^{T}\|A(t,\cdot)\|^{2}_{L^{2}}\frac{dt}{t^{\alpha/\beta}}\lesssim\int_{0}^{T}\|f(t,\cdot)\|^{2}_{L^{2}}\frac{dt}{t^{\alpha/\beta}}.
\end{equation}
\end{lemma}
\begin{lemma}\label{le5}
For $\beta\in(1/2,1)$ and $N(t,x)$ defined on $(0,1)\times
\mathbb{R}^{n},$ let $A(N)$ be the quantity
$$A(\alpha,\beta, N)=\sup_{x\in \mathbb{R}^{n},r\in (0,1)}r^{2\alpha-n+2\beta-2}\int_{0}^{r^{2\beta}}\int_{|y-x|<r}|f(t,x)|\frac{dxdt}{t^{\alpha/\beta}}.$$
Then  for each $k\in \mathbb{N}_{0}:=\mathbb{N}\cup \{0\}$ there
exists a constant $b(k)$ such that the following inequality holds:
\begin{equation}
\int^{1}_{0}\left\|t^{\frac{k}{2}}(-\triangle)^{\frac{k\beta+1}{2}}e^{-\frac{t}{2}(-\triangle)^{\beta}}\int_{0}^{t}N(s,\cdot)ds\right\|_{L^{2}}^{2}
\frac{dt}{t^{\alpha/\beta}}\leq b(k)A(\alpha,\beta,N)
\int^{1}_{0}\int_{\mathbb{R}^{n}}|N(s,x)|\frac{dxds}{s^{\alpha/\beta}}.
\end{equation}
\end{lemma}
\begin{remark}
Similarly when $k=0$, we can prove the following inequality:
\begin{equation}
\int^{1}_{0}\left\|(-\triangle)^{\frac{1}{2}}e^{-t(-\triangle)^{\beta}}\int_{0}^{t}N(s,\cdot)ds\right\|_{L^{2}}^{2}
\frac{dt}{t^{\alpha/\beta}}\lesssim A(\alpha,\beta,N)
\int^{1}_{0}\int_{\mathbb{R}^{n}}|N(s,x)|\frac{dxds}{s^{\alpha/\beta}}.
\end{equation}

\end{remark}
Now we give the main result of this paper.
\begin{theorem}\label{th4}
(Well-posedness)

(i) The subcritical quasi-geostrophic equation (\ref{eq4}) has a
unique small global mild solution in
$(X^{\beta}_{\alpha}(\mathbb{R}^{2}))^{2}$ for all initial data
$\theta_{0}$ with $\nabla\cdot \theta=0$ and
$\|u_{0}\|_{Q^{\beta,-1}_{\alpha}}$ being small.

(ii) For any $T\in(0,\infty)$, there is an $\varepsilon>0$ such that
the quasi-geostrophic equation (\ref{eq4}) has a unique small mild
solution in $(X^{\beta}_{\alpha}(\mathbb{R}^{2}))^{2}$ on
$(0,T)\times\mathbb{R}^{2}$ when the initial data $u_{0}$ satisfies
$\nabla\cdot u_{0}=0$ and $\|u_{0}\|_{(Q^{\beta,-1}_{\alpha,\
T})^{2}}\leq\varepsilon$. In particular, for all $u_{0}\in
\overline{(VQ^{\beta,-1}_{\alpha})^{2}}$ with $\nabla\cdot u_{0}=0$,
there exists a unique small local mild solution in
$(X^{\beta}_{\alpha,T})^{2}$ on $(0,T)\times\mathbb{R}^{2}$.
\end{theorem}
\begin{proof}
By the Picard contraction principle we only need to prove the
bilinear form $B(u,v)$ is bounded on $X^{\beta}_{\alpha}$. We split
the proof into two parts.

{\it Part I:} $\dot{B}^{0,1}_{\infty}(\mathbb{R}^{2})-$boundedness.
The proof of this part has been given in \cite{Marchand Rieusset}.
For completeness, we give the details. We have
\begin{eqnarray*}
&&\|B(u,v)\|_{\dot{B}^{0,1}_{\infty}(\mathbb{R}^{2})}\\
&&\lesssim\int_{0}^{t}\|e^{-(t-s)(-\Delta)^{\beta}}(\partial_{1}(gR_{2}f)-\partial_{2}(gR_{1}f))\|_{\dot{B}^{0,1}_{\infty}(\mathbb{R}^{2})}ds\\
&&\lesssim\int_{0}^{t}\frac{C_{\beta}}{(t-s)^{\frac{1}{2\beta}}s^{1+(1-\frac{1}{\beta})}}s^{1-\frac{1}{2\beta}}\|u\|_{\dot{B}^{0,1}_{\infty}(\mathbb{R}^{2})}
s^{1-\frac{1}{2\beta}}\|v\|_{\dot{B}^{0,1}_{\infty}(\mathbb{R}^{2})}ds\\
&&\lesssim\|u\|_{X^{\beta}_{\alpha}}\|v\|_{X^{\beta}_{\alpha}}\int_{0}^{t}\frac{ds}{(t-s)^{\frac{1}{2\beta}}s^{1+(1-\frac{1}{\beta})}}.
\end{eqnarray*}
Because when $\frac{1}{2}<\beta<1$,
$$\int_{0}^{t/2}\frac{1}{(t-s)^{\frac{1}{2\beta}}s^{1+(1-\frac{1}{\beta})}}ds\lesssim t^{\frac{1}{2\beta}-1}$$
and
$$\int^{t}_{t/2}\frac{1}{(t-s)^{\frac{1}{2\beta}}s^{1+(1-\frac{1}{\beta})}}ds
\lesssim
t^{-2+\frac{1}{\beta}}\int^{t}_{t/2}\frac{1}{(t-s)^{\frac{1}{2\beta}}}ds\lesssim
t^{\frac{1}{2\beta}-1}.$$ Then we can get
$$t^{1-\frac{1}{2\beta}}\|B(u,v)\|_{\dot{B}^{0,1}_{\infty}(\mathbb{R}^{2})}\lesssim\|u\|_{X_{\alpha}^{\beta}}\|v\|_{X^{\beta}_{\alpha}}$$
where in the above estimates we have used for
$f\in\dot{B}^{0,1}_{\infty}(\mathbb{R}^{2})$,
$\|R_{j}f\|_{\dot{B}^{0,1}_{\infty}(\mathbb{R}^{2})}\lesssim\|f\|_{\dot{B}^{0,1}_{\infty}(\mathbb{R}^{2})}$.
In fact by Bernstein's inequality, we have
\begin{eqnarray*}
\sum_{l}\|\Delta_{l}R_{j}f\|_{L^{\infty}(\mathbb{R}^{2})}&=&\sum_{l}\|\partial_{j}(-\Delta)^{-1/2}\Delta_{l}f\|_{L^{\infty}(\mathbb{R}^{2})}\\
&\lesssim&\sum_{l}2^{l}\|(-\Delta)^{-1/2}\Delta_{l}f\|_{L^{\infty}(\mathbb{R}^{2})}\\
&\lesssim&\sum_{l}2^{l}2^{-l}\|\Delta_{l}f\|_{L^{\infty}(\mathbb{R}^{2})}\\
&\leq&\|f\|_{\dot{B}^{0,1}_{\infty}(\mathbb{R}^{2})}.
\end{eqnarray*}
On the other hand, by Young's inequality, we have
\begin{eqnarray*}
&&t^{1-\frac{1}{2\beta}}\|e^{-t(-\Delta)^{\beta}}u_{0}\|_{\dot{B}^{0,1}_{\infty}(\mathbb{R}^{2})}\lesssim
\|u_{0}\|_{\dot{B}^{1-2\beta,\infty}_{\infty}(\mathbb{R}^{2})}\leq\|u_{0}\|_{Q^{\beta,-1}_{\alpha}(\mathbb{R}^{2})}.
\end{eqnarray*}
{\it Part II:} $L^{2}(\mathbb{R}^{2})$-boundedness. Now we estimate
the operation of $B(u,v)$ on the Carleson part of
$X^{\beta}_{\alpha}$. We split again the estimate into two steps.

{\it Step I:} We want to prove the following estimate:
$$r^{2\alpha-2+2\beta-2}\int_{0}^{r^{2\beta}}\int_{|x-y|<r}|B(u,v)|^{2}\frac{dydt}{t^{\alpha/ \beta}}
\lesssim\|u\|_{X^{\beta}_{\alpha}}\|v\|_{X^{\beta}_{\alpha}}.$$ By
symmetry, we only need to deal with the term
\begin{eqnarray*}
&&\int_{0}^{t}e^{-(t-s)(-\Delta)^{\beta}}[\partial_{1}(vR_{1}u)]ds=B_{1}(u,v)+B_{2}(u,v)+B_{3}(u,v)
\end{eqnarray*}
where
$$B_{1}(u,v)=\int_{0}^{t}e^{-(t-s)(-\Delta)^{\beta}}\partial_{1}[(1-1_{r,x})vR_{1}u]ds,$$
$$B_{2}(u,v)=(-\Delta)^{-1/2}\partial_{1}\int_{0}^{t}e^{-(t-s)(-\Delta)^{\beta}}(-\Delta)((-\Delta)^{1/2}(I-e^{-s(-\Delta)^{\beta}})(1_{r,x})vR_{1}u)ds$$
and
$$B_{3}(u,v)=(-\Delta)^{-1/2}\partial_{1}(-\Delta)^{1/2}e^{-t(-\Delta)^{\beta}}\int_{0}^{t}(1_{r,x})vR_{1}uds.$$
For $B_{1}$. Because the n dimensional fractional heat kernel
satisfies the following estimate:
\begin{equation}\label{eq5}
|\nabla
e^{-t(-\Delta)^{\beta}}(x,y)|\lesssim\frac{1}{t^{\frac{n+1}{2\beta}}}\frac{1}{\left(1+\frac{|x-y|}{t^{1/2\beta}}\right)^{n+1}}
\lesssim\frac{1}{(t^{\frac{1}{2\beta}}+|x-y|)^{n+1}},
\end{equation}
we have, for $0<t<r^{2\beta}$ and taking $n=2$ in (\ref{eq5}),
\begin{eqnarray*}
&&|B_{1}(u,v)(t,x)|\\
&&\lesssim\int_{0}^{t}\int_{|z-x|\geq10r}\frac{|R_{1}u(s,z)||v(s,z)|}{|x-z|^{2+1}}dzds\\
&&\lesssim\left(\int_{0}^{r^{2\beta}}\int_{|z-x|\geq10r}\frac{|R_{1}u(s,z)|^{2}}{|x-z|^{3}}dzds\right)^{1/2}
\left(\int_{0}^{r^{2\beta}}\int_{|z-x|\geq10r}\frac{|v(s,z)|^{2}}{|x-z|^{3}}dzds\right)^{1/2}\\
&&:=I_{1}\times I_{2}.
\end{eqnarray*}
For $I_{1}$, we have
\begin{eqnarray*}
I_{1}&\lesssim&\left(\sum_{k=3}^{\infty}\frac{1}{(2^{k}r)^{3}}
\int_{0}^{r^{2\beta}}\int_{|x-z|\leq2^{k+1}r}|R_{1}u(s,x)|^{2}dsdx\right)^{1/2}\\
&\lesssim&\left(\sum_{k=3}^{\infty}\frac{1}{(2^{k}r)^{3}}(2^{k}r)^{2\alpha+2\beta-2}(2^{k}r)^{2-2\beta}
\int_{0}^{r^{2\beta}}\int_{|x-z|\leq2^{k+1}r}|R_{1}u(s,x)|^{2}\frac{dsdx}{s^{\alpha/\beta}}\right)^{1/2}\\
&\lesssim&\left(\sum_{k=3}^{\infty}\frac{(2^{k}r)^{2-2\beta}}{(2^{k}r)}(2^{k}r)^{2\alpha-2+2\beta-2}
\int_{0}^{r^{2\beta}}\int_{|x-z|\leq2^{k+1}r}|R_{1}u(s,x)|^{2}dsdx\right)^{1/2}\\
&\lesssim&\|u\|_{X^{\beta}_{\alpha}}(\sum_{k=3}^{\infty}\frac{1}{2^{k(2\beta-1)}}\frac{1}{r^{2\beta-1}})^{1/2}\\
&\lesssim&\left(\frac{1}{r^{2\beta-1}}\right)^{1/2}\|u\|_{X^{\beta}_{\alpha}}.
\end{eqnarray*}
Similarly we can get
$I_{2}\lesssim\left(\frac{1}{r^{2\beta-1}}\right)^{1/2}\|v\|_{X_{\alpha}^{\beta}}$
 and therefore we have
$$|B_{1}(u,v)|\lesssim\frac{1}{r^{2\beta-1}}\|u\|_{X^{\beta}_{\alpha}}\|v\|_{X^{\beta}_{\alpha}}.$$
Then we can get, using $0<\alpha<\beta$,
\begin{eqnarray*}
\int_{0}^{r^{2\beta}}\int_{|x-y|<r}|B_{1}(u,v)|^{2}\frac{dydt}{t^{\alpha/\beta}}
&\lesssim&
\frac{1}{r^{4\beta-2}}r^{2}\int_{0}^{r^{2\beta}}\frac{dt}{t^{\alpha/\beta}}\|u\|^{2}_{X^{\beta}_{\alpha}}\|v\|^{2}_{X^{\beta}_{\alpha}}\\
&\lesssim&\frac{1}{r^{4\beta-2}}r^{2}r^{2\beta-2\alpha}\|u\|^{2}_{X^{\beta}_{\alpha}}\|v\|^{2}_{X^{\beta}_{\alpha}}\\
&\lesssim&r^{2-2\alpha-2\beta+2}\|u\|^{2}_{X^{\beta}_{\alpha}}\|v\|^{2}_{X^{\beta}_{\alpha}}.
\end{eqnarray*}
That is to say
$$r^{2\alpha-2+2\beta-2}\int_{0}^{r^{2\beta}}\int_{|x-y|<r}|B_{1}(u,v)(t,y)|^{2}\frac{dydt}{t^{\alpha/\beta}}
\lesssim\|u\|^{2}_{X^{\beta}_{\alpha}}\|v\|^{2}_{X^{\beta}_{\alpha}}.$$
For $B_{2}$. By the $L^{2}-$boundedness of Riesz transform, we have
\begin{eqnarray*}
&&\int_{0}^{r^{2\beta}}\int_{|x-y|<r}|B_{2}(u,v)|^{2}\frac{dydt}{t^{\alpha/\beta}}\\
&\lesssim&\int^{r^{\beta}}_{0}
\left\|\int_{0}^{t}e^{-(t-s)(-\Delta)^{\beta}}(-\Delta)((-\Delta)^{-1/2}(I-e^{-s(-\Delta)^{\beta}})(1_{r,x})vR_{1}u)ds\right\|_{L^{2}}^{2}
\frac{dt}{t^{\alpha/\beta}}\\
&\lesssim&\int^{r^{\beta}}_{0}
\left\|\int_{0}^{t}e^{-(t-s)(-\Delta)^{\beta}}(-\Delta)^{\beta}((-\Delta)^{1/2-\beta}(I-e^{-s(-\Delta)^{\beta}})(1_{r,x})vR_{1}u)ds\right\|_{L^{2}}^{2}
\frac{dt}{t^{\alpha/\beta}}\\
&\lesssim&\int_{0}^{r^{2\beta}}t^{2-\frac{1}{\beta}}\int_{|y-x|<r}|R_{1}u(t,y)|^{2}|v(t,y)|^{2}\frac{dydt}{t^{\alpha/\beta}}\\
&\lesssim&\left(\sup_{t>0}t^{1-\frac{1}{2\beta}}\|R_{1}u(t,\cdot)\|_{L^{\infty}(\mathbb{R}^{2})}\right)
\left(\sup_{t>0}t^{1-\frac{1}{2\beta}}\|v(t,\cdot)\|_{L^{\infty}(\mathbb{R}^{2})}\right)\\
&\times&\int_{0}^{r^{2\beta}}\int_{|y-x|<r}|R_{1}u(t,y)||v(t,y)|\frac{dtdy}{t^{\alpha/\beta}}.
\end{eqnarray*}
On one hand, by Bernstein's inequality, we have
$$\|R_{1}u(t,\cdot)\|_{L^{\infty}(\mathbb{R}^{2})}\leq\|R_{1}u(t,\cdot)\|_{\dot{B}^{0,1}_{\infty}(\mathbb{R}^{2})}
\lesssim \|u(t,\cdot)\|_{\dot{B}^{0,1}_{\infty}(\mathbb{R}^{2})}.$$
Then we get
$$\sup_{t>0}t^{1-\frac{1}{2\beta}}\|R_{1}u(t,\cdot)\|_{L^{\infty}(\mathbb{R}^{2})}\lesssim\sup_{t>0}t^{1-\frac{1}{2\beta}}
\|u(t,\cdot)\|_{\dot{B}^{0,1}_{\infty}(\mathbb{R}^{2})}.$$ On the
other hand, we have, by H\"{o}lder's inequality,
\begin{eqnarray*}
&&\int_{0}^{r^{2\beta}}\int_{|x-y|<r}|R_{1}u(t,y)||v(t,y)|\frac{dtdy}{t^{\alpha/\beta}}\\
&\lesssim&\left(\int_{0}^{r^{2\beta}}\int_{|y-x|<r}|R_{1}u(t,y)|^{2}\frac{dtdy}{t^{\alpha/\beta}}\right)^{1/2}
\left(\int_{0}^{r^{2\beta}}\int_{|y-x|<r}|v(t,y)|^{2}\frac{dtdy}{t^{\alpha/\beta}}\right)^{1/2}\\
&\lesssim&r^{2-2\alpha-2\beta+2}\|u\|^{2}_{X_{\alpha}^{\alpha}}\|v\|^{2}_{X_{\alpha}^{\alpha}}.
\end{eqnarray*}
Hence we get
$$\int_{0}^{r^{2\beta}}\int_{|x-y|<r}|B_{2}(u,v)(t,y)|^{2}\frac{dydt}{t^{\alpha/\beta}}
\lesssim
r^{2-2\alpha-2\beta+2}\|u\|^{2}_{X^{\beta}_{\alpha}}\|v\|^{2}_{X^{\beta}_{\alpha}}.$$
For $B_{3}(u,v)$. We have
\begin{eqnarray*}
&&\int_{0}^{r^{2\beta}}\int_{|y-x|<r}|B_{3}(u,v)(t,y)|^{2}\frac{dydt}{t^{\alpha/\beta}}\\
&=&\int_{0}^{r^{2\beta}}\int_{|y-x|<r}\left|(-\Delta)^{-1/2}\partial_{1}(-\Delta)^{1/2}e^{-t(-\Delta)^{\beta}}
\left(\int_{0}^{t}(1_{r,x})vR_{1}udh\right)\right|^{2}\frac{dydt}{t^{\alpha/\beta}}\\
&\lesssim&\int_{0}^{r^{2\beta}}\left\|(-\Delta)^{1/2}e^{-t(-\Delta)^{\beta}}
\left(\int_{0}^{t}(1_{r,x})vR_{1}udh\right)\right\|\frac{dt}{t^{\alpha/\beta}}\\
&\lesssim&r^{2-2\alpha+6\beta-2}\left(\int_{0}^{1}\|M(r^{2\beta}s,\
r\cdot)\|_{L^{1}(\mathbb{R}^{2})}\frac{ds}{s^{\alpha/\beta}}\right)C(\alpha,\beta,f)\\
&\lesssim&r^{2-2\alpha+6\beta-2}r^{2-4\beta}r^{2-4\beta}\|u\|_{X^{\beta}_{\alpha}}\|v\|_{X^{\beta}_{\alpha}}\\
&\lesssim&r^{2-2\alpha-2\beta+2}\|u\|_{X^{\beta}_{\alpha}}\|v\|_{X^{\beta}_{\alpha}}.
\end{eqnarray*}
{\it Step II:} For $j=1,2$, $R_{j}$ are the Riesz transforms
$\partial_{j}(-\Delta)^{-1/2}$. We want to prove:
\begin{equation}\label{eq6}
r^{2\alpha-2+2\beta-2}\int_{0}^{r^{2\beta}}\int_{|x-y|<r}|R_{j}B(u,v)|^{2}\frac{dydt}{t^{\alpha/
\beta}}\lesssim\|u\|_{X^{\beta}_{\alpha}}\|v\|_{X^{\beta}_{\alpha}}.
\end{equation}
 Similar to Step I, we can split $B(u,v)$ into $B_{i}(u,v),\
(i=1,2,3)$. We denote by $A_{i}, i=1,2,3$
\begin{equation}\label{eq7}
A_{i}:=r^{2\alpha-2+2\beta-2}\int_{0}^{r^{2\beta}}\int_{|x-y|<r}|R_{j}B_{i}(u,v)|^{2}\frac{dydt}{t^{\alpha/
\beta}}\lesssim\|u\|_{X^{\beta}_{\alpha}}\|v\|_{X^{\beta}_{\alpha}}.
\end{equation}
In order to estimate the term $A_{1}$, we need the following lemma.
\begin{lemma}\label{le2}
For $\beta>0$, if we denote by $K^{\beta}_{j}$ the kernel of the
operator $e^{-t(-\Delta)^{\beta}}R_{j}$, we have
$$(1+|x|)^{n+|\alpha|}\partial^{\alpha}e^{-t(-\Delta)^{\beta}}R_{j}\in L^{\infty}(\mathbb{R}^{2}).$$
\end{lemma}
\begin{proof}
By Fourier transform, we have
$K^{\beta}_{j}=\mathcal{F}^{-1}(\frac{\xi_{j}}{|\xi|}e^{-|\xi|^{2\beta}})$
where $\mathcal{F}^{-1}$ denotes the inverse Fourier transform.
Because
$$\left[\partial^{\alpha}K_{j}^{\beta}(x)\right]\hat{}\ (\xi)=\frac{\xi_{j}}{|\xi|}|\xi|^{\alpha}e^{-|\xi|^{2\beta}}\in L^{1}(\mathbb{R}^{2}),$$
we have
$$|\partial^{\alpha}K_{j}^{\beta}(x)|\leq\int_{\mathbb{R}^{2}}\left|\frac{\xi_{j}}{|\xi|}|\xi|^{\alpha}e^{-|\xi|^{2\beta}}\right|d\xi\leq C.$$
Then $\partial^{\alpha}K^{\beta}_{j}(x)\in L^{\infty}$. If
$|x|\leq1$, we have
$$(1+|x|)^{n+|\alpha|}|K^{\beta}_{j}(x)|\lesssim C_{\alpha}|K^{\beta}_{j}(x)|\lesssim C.$$
If $|x|>1$, by Littlewood-Paley decomposition and write
$$K^{\beta}_{j}(x)=(Id-S_{0})K^{\beta}_{j}+\sum_{l<0}\Delta_{l}K^{\beta}_{j},$$
where $(Id-S_{0})K^{\beta}_{j}\in \mathcal{S}(\mathbb{R}^{2})$ and
$\Delta_{l}K^{\beta}_{j}=2^{2l}\omega_{j,l}(2^{l}x)$ where
$\widehat{\omega_{j,l}}(\xi)=\psi(\xi)\frac{\xi_{j}}{|\xi|}e^{-|2^{l}\xi|^{2\beta}}\in
L^{1}(\mathbb{R}^{2}).$ Then $\omega_{j,l}(x)_{(l<0)}$ are a bounded
set in $\mathcal{S}(\mathbb{R}^{2}).$ So we have
$$(1+2^{l}|x|)^{N}2^{l(2+|\alpha|)}|\partial^{\alpha}\Delta_{l}K^{\beta}_{j}(x)|\lesssim C_{N}$$
and
\begin{eqnarray*}
&&|\partial^{\alpha}S_{0}K_{j}^{\beta}(x)|\\
&\lesssim&C\sum_{2^{l}|x|\leq1}2^{l(2+|\alpha|)}+\sum_{2^{l}|x|>1}2^{l(2+|\alpha|-N)}|x|^{-N}\\
&\lesssim& C|x|^{-(2+|\alpha|)}.
\end{eqnarray*}
Therefore we have completed the proof of Lemma \ref{le2}\end{proof}
Now we complete the proof of Theorem \ref{th4}. In Lemma \ref{le2},
we take $\alpha=1$ and get
$$\left|\partial_{x}R_{j}e^{-t(-\Delta)^{\beta}}(x,y)\right|\lesssim\frac{1}{(t^{\frac{1}{2\beta}}+|x-y|)^{n+1}}.$$
Similar to the proof in Part I, we can get
$$A_{1}:=r^{2\alpha-2+2\beta-2}\int_{0}^{r^{2\beta}}\int_{|x-y|<r}|R_{j}B_{1}(u,v)|^{2}\frac{dydt}{t^{\alpha/
\beta}}\lesssim\|u\|_{X^{\beta}_{\alpha}}\|v\|_{X^{\beta}_{\alpha}}.$$
In the proof of Theorem \ref{th2}, we in fact prove the following
estimate: for $j=1,2$,
\begin{eqnarray*}
&&r^{2\alpha-2+2\beta-2}\int_{0}^{r^{2\beta}}\int_{|y-x_{0}|<r}\left|R_{j}f(t,y)\right|^{2}\frac{dydt}{t^{\alpha/\beta}}\\
&\lesssim&
\sup_{r>0,x_{0}\in\mathbb{R}^{n}}r^{2\alpha-2+2\beta-2}\int_{0}^{r^{2\beta}}\int_{|y-x_{0}|<r}\left|f(t,y)\right|^{2}\frac{dydt}{t^{\alpha/\beta}}.
\end{eqnarray*}
By the above estimate, we have
\begin{eqnarray*}
&&A_{i}:=r^{2\alpha-2+2\beta-2}\int_{0}^{r^{2\beta}}\int_{|x-y|<r}|R_{j}B_{i}(u,v)|^{2}\frac{dydt}{t^{\alpha/
\beta}}\\
&\lesssim&r^{2\alpha-2+2\beta-2}\int_{0}^{r^{2\beta}}\int_{|x-y|<r}|B_{i}(u,v)|^{2}\frac{dydt}{t^{\alpha/
\beta}}
\end{eqnarray*}
where $i=2,3$. Following the estimate to $B_{i}, i=2,3$, we can get
$$A_{i}:=r^{2\alpha-2+2\beta-2}\int_{0}^{r^{2\beta}}\int_{|x-y|<r}|R_{j}B_{i}(u,v)|^{2}\frac{dydt}{t^{\alpha/
\beta}}\lesssim\|u\|_{X^{\beta}_{\alpha}}\|v\|_{X^{\beta}_{\alpha}}.$$
This completes the proof of Theorem \ref{th4}.

\end{proof}
Following the method applied in  Section 5 of \cite{Li Zhai}, we can
easily get the regularity of the solution to the quasi-geostrophic
equation (\ref{eq1}). So we only state the result and leave the
proof to the readers. For convenience of the study, we introduce a
class of spaces $X^{\beta,k}_{\alpha}$ as follows.
\begin{definition}\label{X k space}
For a nonnegative integer $k$ and $\beta\in(1/2,1],$ we introduce
the space $X^{\beta,k}_{\alpha}$ which is equipped with the
following norm:
$$\|u\|_{X^{\beta,k}_{\alpha}}=\|u\|_{N^{\beta,k}_{\alpha,\infty}}+\|u\|_{N^{\beta,k}_{\alpha,C}}$$
where
\begin{eqnarray*}
\|u\|_{N^{\beta,k}_{\alpha,\infty}}&=&\sup_{\alpha_{1}+\cdots+\alpha_{n}=k}\sup_{t}t^{\frac{2\beta-1+k}{2\beta}}\|\partial_{x_{1}}^{\alpha_{1}}
\cdots\partial_{x_{n}}^{\alpha_{n}}u(\cdot,t)\|_{\dot{B}^{0,1}_{\infty}(\mathbb{R}^{n})},\\
\|u\|_{N^{\beta,k}_{\alpha,C}}&=&\sup_{\alpha_{1}+\cdots+\alpha_{n}=k}\sup_{x_{0},r}\left(r^{2\alpha-n+2\beta-2}\int_{0}^{r^{2\beta}}\int_{|y-x_{0}|<r}
|t^{\frac{k}{2\beta}}\partial_{x_{1}}^{\alpha_{1}}
\cdots\partial_{x_{n}}^{\alpha_{n}}u(t,y)|^{2}\frac{dy
dt}{t^{\alpha/\beta}}\right)^{1/2}\\
&+&\sum^{2}_{j=1}\sup_{\alpha_{1}+\cdots+\alpha_{n}=k}\sup_{x_{0},r}\left(r^{2\alpha-n+2\beta-2}\int_{0}^{r^{2\beta}}\int_{|y-x_{0}|<r}
|R_{j}t^{\frac{k}{2\beta}}\partial_{x_{1}}^{\alpha_{1}}
\cdots\partial_{x_{n}}^{\alpha_{n}}u(t,y)|^{2}\frac{dy
dt}{t^{\alpha/\beta}}\right)^{1/2}.\\
\end{eqnarray*}
\end{definition}
Now we state the regularity result.
\begin{theorem}\label{threg}
 Let $\alpha>0$ and $\max\{\alpha,1/2\}<\beta<1$ with $\alpha+\beta-1\geq0$. There exists an $\varepsilon=\varepsilon(2)$ such that if
$\|u_{0}\|_{Q^{\beta,-1}_{\alpha;\infty}(\mathbb{R}^{2})}<\varepsilon,$
the solution $u$ to equations (\ref{eq1}) verifies:
$$t^{\frac{k}{2\beta}}\nabla^{k}u\in X^{\beta,0}_{\alpha}\text{ for any } k\geq 0.$$
\end{theorem}

\end{document}